%Final version 4 December 2006
%submitted 5 dec 2006
%Typo on page 6 corrected: last display W_x(p) instead of W_x(n) in the submitted version.
%Accepted by the annals 15 October 07;  Revision 8 February 2008.
\input amstex
 \documentstyle{amsppt}
\NoBlackBoxes
\magnification1200
\pagewidth{6.5 true in} 
\pageheight{9.25 true in}
\topmatter
\title
Moments of the Riemann zeta-function
\endtitle 
\author K. Soundararajan 
\endauthor 
\thanks 
The author is partially supported by the National Science Foundation and the American 
Institute of Mathematics (AIM).
\endthanks
\address 
Department of Mathematics, 450 Serra Mall, Bldg. 380, Stanford University, Stanford, CA 94305, USA
\endaddress
\email 
  ksound{\@}math.stanford.edu
\endemail
\dedicatory  
In memoriam Atle Selberg
\enddedicatory
\endtopmatter
\def\lam{\lambda}

\def\Lam{\Lambda}
\def\lam{\lambda}
\document
\loadbold
 \document
 
 \head 1. Introduction \endhead 
 
 \noindent An important problem in analytic number theory is to gain an understanding of the moments
 $$
 M_k(T) = \int_0^T |\zeta(\tfrac 12+it)|^{2k} dt. 
 $$ 
 For positive real numbers $k$, it is believed that 
 $M_k(T) \sim C_k T(\log T)^{k^2}$ for a positive constant $C_k$.  A precise 
 value for $C_k$ was conjectured by Keating and Snaith [10] based on considerations from 
 random matrix theory.  Subsequently, an alternative approach, based on multiple Dirichlet 
 series and producing the same conjecture, was given by Diaconu, Goldfeld and 
 Hoffstein [4].  Recent work by Conrey {\sl et al} [2] gives a more precise conjecture, 
 identifying lower order terms in an asymptotic expansion for $M_k(T)$.   
 
 Despite many 
 attempts, asymptotic formulae for $M_k(T)$ have been established only for $k=1$ 
 (due to Hardy and Littlewood, see [22]) and $k=2$ (due to Ingham, see [22]).  However we do have 
 the lower bound $M_k(T) \gg_k T (\log T)^{k^2}$.  This was 
 established by Ramachandra [13] for positive integers $2k$, by Heath-Brown [6] for all
 positive rational numbers $k$, and assuming the truth of the Riemann Hypothesis 
 by Ramachandra [12] for all positive real numbers $k$.  
 See also the elegant note [1] giving such a bound assuming RH, and 
 [20] for the best known constants implicit in these lower bounds.  
 Analogous conjectures exist (see [2], [4], [11]) for moments of central values of $L$-functions in families, 
 and in many cases  lower bounds of the conjectured order are known (see [16] and [17]).

Here we study the problem of obtaining upper bounds for $M_k(T)$.  When $0\le k \le 2$, Ramachandra ([13], [14]) and  Heath-Brown ([6], [7]) showed, assuming RH,  that $M_k(T) \ll T(\log T)^{k^2}$.  
The Lindel{\" of} Hypothesis is equivalent to the estimate  $M_k(T) \ll_{k,\epsilon} T^{1+\epsilon}$ 
for all natural numbers $k$.   Thus, for $k$ larger than $2$, it seems difficult to make unconditional 
 progress on bounding $M_k(T)$.    
 If we assume RH, then a classical bound of Littlewood (see [22]) gives that (for $t\ge 10$ and some positive 
 constant $C$)
 $$
 |\zeta(\tfrac 12+it)| \ll \exp\Big(C \frac{\log t}{\log \log t}\Big),
\tag{1} $$
 and therefore $M_k(T) \ll T \exp(2kC \log T/\log \log T)$.  We improve 
 upon this, nearly obtaining an upper bound of the conjectured order of magnitude. 
 
 \proclaim{Corollary A}  Assume RH.  For 
 every positive real number $k$, and every $\epsilon >0$ we have 
 $$
T(\log T)^{k^2} \ll_k M_k(T) \ll_{k,\epsilon} T (\log T)^{k^2+\epsilon}.
 $$
 %Moreover, 
 %$$
 %M_k(T) \sim \int_{{\Cal D}} |\zeta(\tfrac 12+it)|^{2k} dt 
 %$$
 \endproclaim 
 
 Our proof suggests that the dominant contribution to the $2k$-th 
 moment comes from $t$ such that $|\zeta(\tfrac 12+it)|$ has 
 size $(\log T)^{k}$, and this set has measure about $T/(\log T)^{k^2}$.

 \proclaim{Corollary B}  Assume RH.  Let $k\ge 0$ be a 
 fixed real number.  For large $T$ we have 
 $$
 \text{\rm meas}\{ t\in [0,T]: |\zeta(\tfrac 12+it)| \ge (\log T)^{k} \} 
 = T (\log T)^{-k^2 +o(1)}.
 $$
 \endproclaim

 We will deduce these corollaries  by finding an upper bound on the 
 frequency with which large values of $|\zeta(\tfrac 12+it)|$ can occur.   
 Throughout we define 
 $$
 {\Cal S}(T,V) = \{ t\in [T,2T]: \  \  \log |\zeta(\tfrac 12+it)| \ge V \}, 
 $$
 and observe that 
 $$ 
 \int_T^{2T} |\zeta(\tfrac 12+it)|^{2k} dt 
 = -\int_{-\infty}^{\infty} e^{2kV} d \text{meas}({\Cal S}(T,V))  
 = 2k\int_{-\infty}^{\infty} e^{2kV} \text{meas}({\Cal S}(T,V)) dV. \tag{2} 
 $$
To prove Corollaries A and B, we desire estimates for the measure of ${\Cal S}(T,V)$ 
 for large $T$ and all $V\ge 3$. 
 To place our main theorem in context, let us recall the beautiful result 
 of Selberg that as $t$ varies in $[T,2T]$ the distribution of 
 $\log |\zeta(\tfrac 12+it)|$ is approximately Gaussian with mean $0$ and 
 variance $\tfrac 12 \log \log T$.  Precisely, Selberg's theorem (see [18], [19]) shows that for 
 any fixed $\lambda \in {\Bbb R}$, and as $T\to \infty$ 
 $$ 
\text{meas} ({\Cal S}(T,\lambda \sqrt{\tfrac 12 \log \log T})) = T\Big(\frac{1}{\sqrt{2\pi} }\int_{\lambda}^{\infty} e^{-x^2/2} 
 dx + o(1)\Big).
 $$
 Although Selberg's result holds only for $V$ of size $\sqrt{\log \log T}$, we may speculate 
 that in a much larger range for $V$ a similar estimate holds:
 $$ 
\text{meas}( {\Cal S}(T,V)) \ll T \frac{\sqrt{\log \log T}}{V} \exp\Big( -\frac{V^2}{\log \log T}\Big).  \tag{3} 
 $$
 Such an estimate would lead, via (2), to the bound $M_k(T) \ll_k T(\log T)^{k^2}$.  In our 
 main Theorem we establish that a weaker form of (3) holds in the 
 range $3 \le V = o((\log \log T )\log_3 T)$, where throughout 
 we shall write $\log_3$ for $\log \log \log$.   In the application to moments, the crucial range 
 is when $V$ is of size about $k\log \log T$, and our Theorem shows that a version of 
 (3) holds for such $V$.   For larger values of $V$ we obtain 
 an upper bound of the form $\text{meas}({\Cal S}(T,V)) \ll T\exp(-\frac {1}{33} V\log V)$ 
 (at least when $\frac 12 (\log \log T) \log_3 T <V$), and the shape of this estimate is in keeping with 
 Littlewood's bound (1).

 \proclaim{Theorem} Assume RH.  Let $T$ be large, and $V\ge 3$ be a real number, 
 and let ${\Cal S}(T,V)$ be as defined above.  %Consider the set 
% $$
%{\Cal S}={\Cal S}(T,V)= \{ t\in [T,2T]: \ \ \log |\zeta(\tfrac 12+it)| \ge V\}. 
 %$$
If $10\sqrt{\log \log T} \le V\le \log \log T $ then 
we have 
$$
\text{\rm meas}({\Cal S}(T,V)) \ll T \frac{V}{\sqrt{\log \log T} } \exp\Big( -\frac{V^2}{\log \log T} \Big( 1- \frac{4}{ \log_3 T}\Big)\Big);
$$ 
 if $\log \log T < V \le \frac 12 (\log \log T) \log_3 T$ we have 
$$ 
\text{\rm meas}({\Cal S}(T,V)) \ll T \frac{V}{\sqrt{\log \log T}} \exp\Big( - \frac{V^2}{\log \log T} 
  \Big( 1- \frac{7V}{4(\log \log T) \log_3 T}\Big)^{2} \Big); 
$$
and, finally, if $\frac 12(\log \log T) \log_3 T < V$ we have 
$$ 
\text{\rm meas}({\Cal S}(T,V)) \ll T \exp\Big( -\frac{1}{33} V\log V\Big). 
$$ 
%  The measure of the set ${\Cal S}(T,V)$ is bounded by 
 %$$
 %\ll T \Big( \exp \Big( -\frac{V^2}{\log \log T} \Big) + \exp \Big(-\frac 43 V\log V \Big)\Big). 
 %$$
 % following bounds for the measure of ${\Cal S}$ hold: 
% when $V\le \log \log T$ we have
 %$$
 %\text{\rm {meas}}({\Cal S}) \ll T \exp\Big(-\frac{V^2}{\log \log T} \Big(1-\frac{300}{\log_3 T}\Big)\Big),
% $$
% when $\log \log T \le V \le (\log \log T \log_3 T)/1200$ we have 
 %$$
 %\text{\rm meas}({\Cal S}) \ll T \exp\Big(-\frac{V^2}{\log \log T} 
 %\Big(1-\frac{300V}{\log \log T\log_3 T}\Big)\Big),
% $$
% while if $V>(\log \log T \log_3 T)/1200$ we have 
 %$$
 %\text{\rm meas}({\Cal S}) \ll T\exp\Big(-\frac{V\log V}{1000}\Big).
 %$$
 \endproclaim

 In the limited range $0\le V\le \log\log T$, Jutila [8] had previously shown that 
 $$
 \text{meas}({\Cal S}(T,V)) \ll T\exp\Big(-\frac{V^2}{\log \log T} \Big(1+ O\Big(\frac{V}{\log \log T}\Big)\Big).
 $$ 
 I have shown recently in [21] that in the range 
 $3 \le V\le \frac 15\sqrt{\log T/\log \log T}$ one has 
$$
\text{meas}({\Cal S}(T,V)) \gg \frac{T}{(\log T)^4} \exp\Big(-10\frac{V^2}{\log \frac{\log T}{8V^2 \log V}}\Big).
$$
In contrast to our Theorem above, these results are unconditional.

As noted above, we build on Selberg's work on the distribution of $\log \zeta(\tfrac 12+it)$.    
Selberg computed the moments of the real and imaginary parts of $\log \zeta(\tfrac 12+it)$.  
To achieve this he found an ingenious expression for $\log \zeta(\tfrac 12+it)$ in terms of 
primes.  His ideas work very well 
 for the imaginary part of the logarithm, but are more complicated for the real part 
 of the logarithm owing to zeros of the zeta-function lying very close to $\tfrac 12+it$.  
 One novelty in our work is the realization that if we seek an upper bound for $\log |\zeta(\tfrac 12+it)|$ (which is what is needed for our Theorem) then the effect of the zeros very near $\tfrac 12 +it$ is 
 actually benign.   This is our main Proposition given below, from which we will deduce our Theorem.  We should comment that Selberg's work is unconditional, and uses zero-density 
 results which put most of the zeros near the critical line.  It would be interesting to see 
 how much of our work can be recovered unconditionally.

\proclaim{Proposition}  Assume RH.  Let $T$ be large, let $t\in [T,2T]$, and let $2\le x\le T^2$.  Let $\lambda_0=0.4912\ldots$ 
denote the unique positive real number satisfying $e^{-\lambda_0} = \lambda_0+ \lambda_0^2/2$.   
For all $\lambda \ge \lambda_0$ we have the estimate 
$$
\log |\zeta(\tfrac 12+it)| \le \text{\rm Re } \sum_{n\le x} \frac{\Lam(n)}{n^{\frac 12+ \frac{\lam}{\log x} +it} \log n} 
\frac{\log (x/n)}{\log x} + \frac{(1+\lam)}{2} \frac{\log T}{\log x} + O\Big( \frac{1}{\log x}\Big).
$$
\endproclaim 

Taking $x=(\log T)^{2-\epsilon}$ in our Proposition, and estimating the sum over $n$ trivially, 
we obtain the following explicit form of Littlewood's bound (1), which improves upon the 
previous estimate obtained by Ramachandra and Sankaranarayanan [15].   There 
is certainly some scope to improve our Corollary C below, and it may be instructive to understand 
what the limit of the method would be (analogously to the elegant treatment of $\text{Im }\log \zeta(\tfrac 12 +it)$ given by Goldston and Gonek [5]).   
 
 \proclaim{Corollary C}  Assume RH.  For all large $t$ we have 
 $$
 |\zeta(\tfrac 12+it)| \le \exp\Big( \Big(\frac{1+\lambda_0}{4}+o(1)\Big)\frac{\log t}{\log \log t} \Big) 
 \le \exp\Big (\frac{3}{8} \frac{\log t}{\log \log t}\Big).
 $$
 \endproclaim

 The method developed here is robust and applies equally well to moments in families of $L$-functions; we  discuss this briefly in section 4 below.   We end the introduction by deriving Corollaries A and B.

 \demo{Proof of  Corollaries A and B}  As mentioned earlier, the lower bound for $M_k(T)$ 
 in Corollary A is due to Ramachandra [12].   The upper bound follows upon inserting the 
 bounds of the Theorem into (2).  In performing this computation, it is convenient to use 
 our Theorem in the crude form $\text {meas}({\Cal S}(T,V)) \ll T (\log T)^{o(1)} \exp(-V^2/\log \log T)$ 
 for $3 \le V \le 4k \log \log T$, and $\text{meas}({\Cal S}(T,V)) \ll T (\log T)^{o(1)} \exp(-4kV)$ 
 for $V> 4k \log \log T$.  
 
 From the Theorem, we may see that the contribution to $M_k(T)$ 
 from $t\in [0,T]$ with $|\zeta(\tfrac 12+it)|>(\log T)^{k+\epsilon}$ or 
 $|\zeta(\tfrac 12+it)| < (\log T)^{k-\epsilon}$ is $o(T(\log T)^{k^2})$.   Combining this 
 with the lower bound $M_k(T)\gg_k T(\log T)^{k^2}$ we obtain the lower 
 bound for the measure implicit in Corollary B.  The upper bound implicit there follows 
 from the Theorem.  
   \enddemo 
  
  {\bf Acknowledgments.}   I am grateful to Mike Rubinstein for some helpful remarks, and 
  to Brian Conrey for his encouragement. 
  
 \head 2. Proof of the main Proposition \endhead 
 
 \noindent  In proving the Proposition we may suppose that $t$ does not coincide 
 with the ordinate of a zero of $\zeta(s)$.   Letting $\rho=\tfrac 12+i\gamma$ run over the non-trivial 
 zeros of $\zeta(s)$, we define  
 $$
 F(s) = \text{Re }\sum_{\rho} \frac{1}{s-\rho} = \sum_{\rho} \frac{\sigma-1/2}{(\sigma-1/2)^2+ 
 (t-\gamma)^2}. 
 $$
 Visibly $F(s)$ is non-negative in the half-plane $\sigma \ge 1/2$.   Recall 
 Hadamard's factorization formula which gives (see (8) and (11) of \S 12 of Davenport [3]) 
  $$ 
 \text{Re }\frac{\zeta^{\prime}}{\zeta}(s) =   -\text{Re }\frac{1}{s-1} + \frac 12\log \pi - \frac 12 
 \text{Re }\frac{\Gamma^{\prime}}{\Gamma}(\tfrac 12 s+1) + F(s),
 $$
 so that for $t\in [T,2T]$ an application of Stirling's formula yields 
 $$
-\text{Re }\frac{\zeta^{\prime}}{\zeta}(s)  =\frac 12 \log T +O(1) - F(s). \tag{4}
$$
Integrating (4) as $\sigma = \text{Re}(s)$ varies from $\tfrac 12$ to $\sigma_0 (> \frac 12)$ 
we obtain, setting $s_0 =\sigma_0 +it$,  
 $$
\align
 \log |\zeta(\tfrac 12+it)| - \log |\zeta(s_0)| 
& = \Big(\frac {\log T}{2} +O(1)\Big) (\sigma_0-\tfrac 12) -\int_{1/2}^{\sigma_0} F(\sigma+it) d\sigma \\
&=  (\sigma_0-\tfrac 12) \Big( \frac{\log T}{2} + O(1)\Big) - 
\frac{1}{2}\sum_{\rho} \log \frac{(\sigma_0-\frac 12)^2+(t-\gamma)^2}{(t-\gamma)^2}.\\
\endalign
 $$ 
 Since $\log (1+x^2) \ge x^2/(1+x^2)$ we deduce that 
 $$
 \log |\zeta(\tfrac 12+it)| - \log |\zeta(s_0)| 
 \le (\sigma_0-\tfrac 12) \Big(\frac{\log T}{2}+ O(1)-\frac{F(s_0)}{2}\Big). \tag{5}
 $$
 
 \proclaim{Lemma 1}  Unconditionally, for any $s$ not coinciding with $1$ or a 
 zero of $\zeta(s)$ and any $x\ge 2$, we have 
 $$
 \align
 -\frac{\zeta^{\prime}}{\zeta}(s)= \sum_{n\le x} 
 \frac{\Lambda(n)}{n^s} \frac{\log (x/n)}{\log x} &+ \frac{1}{\log x} \Big(\frac{\zeta^{\prime}}{\zeta}(s)\Big)^{\prime} 
 + \frac{1}{\log x} \sum_{\rho} \frac{x^{\rho-s}}{(\rho-s)^2} \\
 &-\frac{x^{1-s}}{(1-s)^2 \log x} 
 + \frac{1}{\log x} \sum_{k=1}^{\infty} \frac{x^{-2k-s}}{(2k+s)^2}.
 \\
 \endalign
 $$
 \endproclaim 
 \demo{Proof}  This is similar to an identity of Selberg, see Theorem 14.20 of Titchmarsh [22].    
 With $c= \max(1,2-\sigma)$ we consider 
 $$
 \frac{1}{2\pi i} \int_{c-i\infty}^{c+i\infty} -\frac{\zeta^{\prime}}{\zeta}(s+w) 
 \frac{x^w}{w^2} dw  = \sum_{n\le x}
  \frac{\Lambda(n)}{n^s} \log (x/n),
 $$  
 upon integrating term by term using the Dirichlet series expansion of 
 $-\frac{\zeta^{\prime}}{\zeta}(s+w)$. 
  On the other hand, moving the line of integration to the left and calculating residues 
  %(from the poles at $w=0$, $w=\rho-s$, $w=1-s$ and $w=-2k-s$ for $k\in {\Bbb N}$)
   this 
 equals 
 $$
 -\frac{\zeta^{\prime}}{\zeta}(s) \log x - \Big(\frac{\zeta^{\prime}}{\zeta}(s)\Big)^{\prime} 
 -\sum_{\rho} \frac{x^{\rho-s}}{(\rho-s)^2} +\frac{x^{1-s}}{(1-s)^2}-\sum_{k=1}^{\infty} 
 \frac{x^{-2k-s}}{(2k+s)^2}. 
 $$
 Equating these two expressions we obtain the lemma.  
 \enddemo
 
Take $s=\sigma+it$ in Lemma 1, extract the real parts of both sides, and integrate 
over $\sigma$ from $\sigma_0$ to $\infty$.  
Thus, for $2\le x\le T^2$, 
 $$
 \align
\log |\zeta(s_0)| = \text{Re }\Big( \sum_{n\le x} \frac{\Lambda(n)}{n^{s_0} \log n} \frac{\log (x/n)}{\log x}  
  &- \frac{1}{\log x} \frac{\zeta^{\prime}}{\zeta}(s_0)\\
  & + \frac{1}{\log x} \sum_{\rho} 
 \int_{\sigma_0}^{\infty} \frac{x^{\rho-s}}{(\rho-s)^2} d\sigma +O\Big(\frac{1}{\log x}\Big)\Big).
 \\
 \endalign
 $$
Using (4), and observing that 
$$
\sum_{\rho}\Big|\int_{\sigma_0}^{\infty} \frac{x^{\rho -s}}{(\rho -s)^2} d\sigma\Big| 
\le \sum_{\rho}\int_{\sigma_0}^{\infty}\frac{ x^{\frac 12-\sigma}}{|s_0-\rho|^2} d\sigma 
= \sum_{\rho}\frac{x^{\frac 12-\sigma_0}}{|s_0-\rho|^2 \log x}= \frac{x^{\frac 12-\sigma_0}F(s_0)}{(\sigma_0-\frac 12)\log x},
$$ 
we deduce that 
 $$
 \align
 \log |\zeta(s_0)| 
 &\le  \text{Re }\sum_{n\le x} \frac{\Lambda(n)}{n^{s_0} \log n} \frac{\log (x/n)}{\log x} 
 + \frac{\log T}{2\log x} -\frac{F(s_0)}{\log x} + \frac{x^{\frac 12-\sigma_0}F(s_0)}{(\sigma_0-1/2)\log^2 x} 
 + O\Big(\frac{1}{\log x}\Big).\tag{6}\\
 \endalign
 $$
 
Adding the inequalities (5) and (6), we obtain that 
 $$
 \align
 \log |\zeta(\tfrac 12+ it)|  &\le 
 \text{Re }\sum_{n\le x} \frac{\Lambda(n)}{n^{s_0} \log n} \frac{\log x/n}{\log x} 
 + \frac{\log T}{2} \Big(\sigma_0 -\frac 12 + \frac{1}{\log x}\Big) 
 \\
 &\hskip .5 in +F(s_0) \Big( \frac{x^{\frac 12-\sigma_0}}{(\sigma_0-\frac 12) \log^2 x} -\frac{1}{\log x} 
 -\frac{(\sigma_0-\frac 12)}{2}\Big) + O\Big(\frac{1}{\log x}\Big).\tag{7} \\
 \endalign
 $$
 We choose $\sigma_0 =\frac12 + \frac{\lam}{\log x}$, where $\lam \ge \lam_0$.  This restriction 
 on $\lam$ ensures that the term involving $F(s_0)$ in (7) makes a negative contribution, 
 and may therefore be omitted.  The Proposition follows.

\head 3. Proof of the Theorem \endhead 

\noindent Our proof of the Theorem rests upon our main Proposition.  We begin 
by showing that the sum over prime powers appearing there may be restricted just 
to primes. 

\proclaim{Lemma 2}  Assume RH.  Let $T\le t\le 2T$, let $2\le x\le T^2$, and let $\sigma\ge \frac 12$.    Then 
$$
\Big| \sum\Sb n\le x \\ n \neq p\endSb \frac{\Lambda(n)}{n^{\sigma+it}\log n}  
\frac{\log x/n}{\log x} \Big| \ll \log \log\log T +O(1).
$$
\endproclaim 
\demo{Proof}   The terms when $n=p^k$ for $k\ge 3$ clearly contribute an amount $\ll 1$, and 
it remains to handle the terms $n=p^2$.  
By following closely the explicit formula proof of the prime number theorem (see \S 17 and 18 
of Davenport [3]) we obtain that, on RH, $\sum_{p\le z} (\log p)p^{-2it}  \ll z/T + \sqrt{z} (\log zT)^2$.
By partial summation, using this estimate when $z\ge (\log T)^4$ and the trivial $\ll z$ 
for smaller $z$, we deduce that for $\sigma \ge \frac 12$ 
$$
\sum_{p\le \sqrt{x}} \frac{1}{p^{2\sigma+2it}} \frac{\log (\sqrt{x}/p)}{\log \sqrt{x}} \ll \log \log \log T, 
$$
completing our proof.  

\enddemo  
 
We also need  a standard mean value estimate whose proof we include 
 for completeness.  
 
 \proclaim{Lemma 3}  Let $T$ be large, and let $2\le x\le T$.  Let $k$ be a natural 
 number such that $x^k\le T/\log T$.   For any complex numbers $a(p)$ we have 
 $$
 \int_T^{2T} \Big| \sum_{p\le x} \frac{a(p)}{p^{\frac 12+it}}\Big|^{2k} dt 
 \ll T k! \Big( \sum_{p\le x} \frac{|a(p)|^2}{p}\Big)^k.
 $$
 \endproclaim 
\demo{Proof}  Write 
$$
\Big( \sum_{p\le x} \frac{a(p)}{p^{\frac 12+it}}\Big)^k = \sum_{n\le x^{k}} \frac{a_{k,x}(n)}{n^{\frac 12+it}},
$$
where $a_{k,x}(n)=0$ unless $n$ is the product of $k$ (not necessarily distinct) primes, all below $x$.  
In that case, if we write the prime factorization of $n$ as $n=\prod_{i=1}^{r} p_i^{\alpha_i}$ then 
$a_{k,x}(n) = \binom{k}{\alpha_1,\ldots,\alpha_r} \prod_{i=1}^{r} a(p_i)^{\alpha_i}$.   Now 
$$
\align
\int_{T}^{2T} \Big| \sum_{p\le x} \frac{a(p)}{p^{\frac 12+it}}\Big|^{2k} dt 
&= \sum_{m,n\le x^k} \frac{a_{k,x}(m)\overline{a_{k,x}(n)}}{\sqrt{mn}} \int_T^{2T} \Big(\frac{n}{m}\Big)^{it} 
dt 
\\
&= T \sum_{n\le x^k} \frac{|a_{k,x}(n)|^2}{n} + O\Big( \sum\Sb m,n\le x^k \\ m\neq n \endSb 
\frac{|a_{k,x}(m){a_{k,x}(n)}|}{\sqrt{mn}|\log (m/n)|}\Big),
\\
\endalign
$$
upon separating the diagonal terms $m=n$, and the off-diagonal terms $m\neq n$.    
Since $2|a_{k,x}(m)a_{k,x}(n)/\sqrt{mn}| \le |a_{k,x}(m)|^2/m+ |a_{k,x}(n)|^2/n$ we see that the 
off-diagonal terms above contribute 
$$
\ll \sum_{n\le x^k }\frac{|a_{k,x}(n)|^2}{n} \sum\Sb m\le x^k\\ m\neq n \endSb \frac{1}{|\log (m/n)|} 
\ll x^k \log (x^k) \sum_{n\le x^k} \frac{|a_{k,x}(n)|^2}{n} \ll T \sum_{n\le x^k} \frac{|a_{k,x}(n)|^2}{n},
$$
recalling that $x^k \le T/\log T$.  The Lemma follows upon noting that
$$
\align
\sum_{n\le x^k} \frac{|a_{k,x}(n)|^2}{n} &= \sum\Sb p_1< \ldots < p_r \le x  \endSb 
\sum\Sb \alpha_1, \ldots, \alpha_r \ge 1\\ \sum\alpha_i =k \endSb 
\binom{k}{\alpha_1,\ldots,\alpha_r}^2 \frac{|a(p_1)|^{2\alpha_1}\cdots |a(p_r)|^{2\alpha_r}}{p_1^{\alpha_1} \cdots p_r^{\alpha_r}} 
\\
&\le k! \sum\Sb p_1< \ldots < p_r \le x \endSb 
\sum\Sb \alpha_1, \ldots, \alpha_r \ge 1\\ \sum\alpha_i =k \endSb 
\binom{k}{\alpha_1,\ldots,\alpha_r} \frac{|a(p_1)|^{2\alpha_1}\cdots |a(p_r)|^{2\alpha_r}}{p_1^{\alpha_1} \cdots p_r^{\alpha_r}} 
\\
&=k! \Big(\sum_{p\le x} \frac{|a(p)|^2}{p}\Big)^{k}.
\\
\endalign
$$

\enddemo 

In proving our Theorem we may assume that $10\sqrt{\log \log T} \le V\le \frac 38 \log T/\log \log T$.  
We also keep in mind that $T$ is large.   We define a parameter $A$ by setting $A=\frac 12 \log_3 T$ 
when $V\le \log \log T$, setting $A= \frac{\log\log T}{2V} \log_3 T$ when $\log \log T < V \le 
\frac 12(\log\log T)\log_3 T$, and finally setting $A=1$ when $V>\frac 12(\log \log T)\log_3 T$.   We 
further set $x= T^{A/V}$ and $z=x^{1/\log \log T}$.   Using Lemma 2 and our Proposition, we 
find that 
$$
\log |\zeta(\tfrac 12+it)| \le S_1(t) + S_2 (t)+ \frac{1+\lambda_0}{2A} V + O(\log \log \log T), %\tag{7a}
$$ 
where 
$$
S_1(t) = \Big| \sum_{p\le z} \frac{1}{p^{\frac 12+\frac{\lam_0}{\log x} +it}}\frac{\log (x/p)}{\log x}\Big|,  %\tag{7b}
$$
and 
$$
S_2(t)= \Big| \sum_{z< p\le x} \frac{1}{p^{\frac 12+\frac{\lam_0}{\log x} +it}}\frac{\log (x/p)}{\log x}\Big|. %\tag{7c}
$$
If $t\in {\Cal S}(T,V)$ then we must either have %$S_2(t) \ge V/(8A)$, or $S_1(t) \ge V(1-1/(8A))$.
$$
S_2(t) \ge \frac{V}{8 A}, \qquad \text{ or } \qquad 
S_1(t) \ge V \Big(1- \frac{7}{8A}\Big) =: V_1,
$$
say.

By Lemma 3 we see that for any natural number $k \le V/A - 1$ we have 
$$
\int_T^{2T} |S_2(t)|^{2k} dt\ll T k! \Big( \sum_{z<p\le x} \frac{1}{p} \Big)^{k} 
\ll T \Big(k (\log \log \log T +O(1))\Big)^k.
$$
Choosing $k$ to be the largest integer below $V/A-1$, we obtain that 
the measure of $t\in [T,2T]$ with $S_2(t) \ge V/(8A)$ is  
$$
\ll T \Big( \frac{8 A}{V}\Big)^{2k} (2k \log \log \log T)^k \ll T \exp\Big( -  \frac{V}{2A}\log V \Big). 
\tag{8}
$$

  Next we consider the measure of the set of $t\in [T, 2T]$ with $S_1(t)\ge V_1$.  
  By Lemma 3, we see that for any natural number $k\le \log (T/\log T)/\log z$, 
  $$
  \int_T^{2T} |S_1(t)|^{2k} dt \ll T k! \Big(\sum_{p\le z} \frac{1}{p} \Big)^{k} 
  \ll T \sqrt{k} \Big( \frac{k\log \log T}{e}\Big)^{k}, 
  $$
  so that the measure of $t \in [T,2T]$ with $S_1(t) \ge V_1$ is 
  $$
  \ll T\sqrt{k} \Big( \frac{k\log \log T}{eV_1^2}\Big)^{k}.
  $$
  When $V\le (\log \log T)^2$ we choose $k$ to be 
   the greatest integer below $V_1^2/\log \log T$, and 
  when $V>(\log \log T)^2$ we choose $k = [10V]$.  
  It then follows that the measure of $t\in [T,2T]$ with 
  $S_1(t) \ge V_1$ is 
  $$
  \ll T \frac{V}{\sqrt{\log \log T}} \exp\Big( -\frac{V_1^2}{\log \log T} \Big) + 
  T\exp ( - 4V\log V). \tag{9}
  $$ 
   
  Our Theorem follows upon combining the estimates (8) and (9).

\head 4.  Moments of $L$-functions in families \endhead
\def\sumstar{\sideset\and^*\to \sum}
\def\sumflat{\sideset\and^\flat \to \sum}

\noindent We briefly sketch here the modifications needed to obtain 
bounds for moments of $L$-functions in families.  Throughout we assume 
the Generalized Riemann Hypothesis for the appropriate $L$-functions 
under consideration.  If $q$ is a large prime, then Rudnick and 
Soundararajan [16] showed that, for positive rational numbers $k$,
$$
\sumstar_{\chi\pmod q} |L(\tfrac 12, \chi)|^{2k} \gg_k q(\log q)^{k^2},
$$
where the sum is over primitive characters $\chi$.  The argument given 
here carries over to obtain the upper 
bound $\ll_{k,\epsilon} q(\log q)^{k^2 +\epsilon}$ for all positive real $k$.  The only difference 
is that one uses the orthogonality relations of the characters $\pmod q$ 
to treat non-diagonal terms in the analog of Lemma 3.  

More interesting is the case of quadratic Dirichlet $L$-functions.  
In [17] Rudnick and Soundararajan showed that for rational numbers $k\ge 1$
$$
\sumflat_{|d|\le X} L(\tfrac 12,\chi_d)^k \gg_k X(\log X)^{k(k+1)/2},
$$
where the sum is over fundamental discriminants $d$, and $\chi_d$ denotes 
the associated primitive quadratic character.  To obtain an upper bound, 
we seek to bound the frequency of large values of $L(\tfrac 12,\chi_d)$.  
Analogously to our Proposition we find that\footnote{Since we are assuming GRH we know that $L(\frac 12,\chi_d)\ge 0$.  If $L(\frac 12,\chi_d)=0$ we interpret $\log L(\frac 12,\chi_d)$ as $-\infty$ so 
that the claimed inequality is still valid.}, for any $x\ge 2$ and with $\lambda_0$ as in 
our Proposition, 
$$
\log L(\tfrac 12,\chi_d) \le \sum_{2\le n\le x} \frac{\Lambda(n)\chi_d(n)}{n^{\frac 12+\frac{\lambda_0}{\log x}}\log n}\frac{\log (x/n)}{\log x} 
+ \frac{(1+\lambda_0)}{2} \frac{\log |d|}{\log x} + O\Big(\frac{1}{\log x}\Big).
$$
Notice that, in contrast to Lemma 2, the contribution of the prime squares in our sum is $\sim \frac 12 \log \log x$, since
$\chi_d(p^2)=1$ for all $p\nmid d$.  Taking this key difference into account, we may argue
as in \S 4, using now quadratic reciprocity and the P{\' o}lya-Vinogradov 
inequality to develop the analog of Lemma 3. Thus we obtain that the number of $d$ with $|d|\le X$ 
and $\log L(\frac 12,\chi_d) \ge V+\frac 12\log \log X$ is 
$$
\ll X\exp\Big(- \frac{V^2}{2\log \log X}(1+o(1))\Big) 
$$ 
when $\sqrt{\log \log X} \le V = o((\log \log X)\log_3 X)$, and when $V\ge (\log \log X)\log_3 X$ this 
number is $\ll X \exp(-cV\log V)$ for some positive constant $c$. From these estimates we deduce that 
$$
\sumflat_{|d|\le X} L(\tfrac 12, \chi_d)^k \ll_{k,\epsilon} X(\log X)^{k(k+1)/2+\epsilon}.
$$

As a last example consider the family of quadratic twists of a given elliptic 
curve $E$.  We write the $L$-function for $E$ as $L(E,s)=\sum_{n=1}^{\infty} a(n)n^{-s}$ 
where the $a(n)$ are normalized to satisfy $|a(n)|\le d(n)$ (the number of divisors of $n$). 
The methods of [16] and [17] can be used to show that for rational 
numbers $k\ge 1$ 
$$
\sumflat_{|d|\le X} L(E\otimes\chi_d, \tfrac 12)^k \gg_{k} X(\log X)^{k(k-1)/2}.
$$
Writing $a(p)=\alpha_p+\beta_p$ with $\alpha_p\beta_p=1$ we obtain analogously to our Proposition  
that, for $x\ge 2$ and $\lambda_0$ as before, 
$$
\align
\log L(E\otimes \chi_d, \tfrac 12) &\le \sum\Sb n=p^{\ell} \le x \\ \ell \ge 1\endSb 
\frac{\chi_d(p^\ell) (\alpha_p^\ell +\beta_p^\ell)}{\ell n^{\frac 12+\frac{\lambda_0}{\log x}}} \frac{\log (x/p^{\ell})}{\log x}  
+ (1+\lambda_0) \frac{\log |d|}{\log x}  + O\Big(\frac{1}{\log x}\Big).
\endalign
$$
The contribution of the terms $\ell \ge 3$ above is $O(1)$, and the contribution of $\ell=2$ 
(the prime squares) is 
$$
\sum\Sb p\le \sqrt{x} \\ p\nmid d\endSb  \frac{a(p^2)-1}{2p^{1+2\frac{\lambda_0}{\log x}}} \frac{\log (x/p^2)}{\log x} 
 \sim -\frac 12 \log \log x.
$$
After taking this feature into account, we may develop the analogous argument of \S 4.  
Thus, the number of $d$ with $|d|\le X$ and $\log L(E\otimes \chi_d,\tfrac 12) \ge V-\frac
12\log \log X$ is 
$$
\ll X \exp\Big( -\frac{V^2}{2\log \log X} (1+o(1)) \Big) 
$$
when $\sqrt{\log \log X} \le V =o((\log \log X)\log_3 X)$ and when $V\ge (\log \log X)\log_3 X$ 
this number is $\ll X\exp(-cV\log V)$ 
for some positive constant $c$.  This leads to the upper bound 
$$
\sumflat_{|d|\le X} L(E\otimes \chi_d,\tfrac 12)^k \ll X(\log X)^{k(k-1)/2+\epsilon}.
$$

Our work above is in keeping with conjectures of Keating and Snaith [11]  (see (51) and (79) there)
that an analog of Selberg's result holds in families of $L$-functions.  
Thus in the unitary family of $\chi \pmod q$, we expect that 
the distribution of $\log |L(\tfrac 12,\chi)|$ is Gaussian with mean $0$, and 
variance $\sim \tfrac 12 \log \log q$.   In the symplectic family of quadratic Dirichlet 
characters, we expect that the distribution of $\log L(\tfrac 12,\chi_d)$ 
is Gaussian with mean $\sim \tfrac{1}{2} \log \log |d|$ and variance 
$\sim \log \log |d|$.  Thus most values of $L(\tfrac 12,\chi_d)$ are quite 
large.  In the orthogonal family of quadratic twists of an elliptic curve, first 
we must restrict to those twists with positive sign of the functional 
equation, else the $L$-value is $0$.  In this restricted class, we expect 
that the distribution of $\log L(E\otimes \chi_d, \tfrac 12)$  is Gaussian with 
mean $\sim -\tfrac 12 \log \log |d|$ and variance $\sim \log \log |d|$.  Thus 
the values in an orthogonal family tend to be small.  With a little more work, the 
ideas in this paper would show (assuming GRH) that the frequency with 
which $\log |L|$ exceeds $\text{Mean} + \lambda \cdot \text{Var}$ is 
bounded above by $\frac{1}{\sqrt{2\pi}}\int_{\lambda}^{\infty} e^{-x^2/2} dx $ for any fixed 
real number $\lambda$.  
If in addition to GRH we assume that most of the $L$-functions under consideration
do not have a zero very near\footnote{This is implied by the one level density conjectures of Katz 
and Sarnak [9].} 
 $\tfrac 12$, then Selberg's techniques would 
yield these Keating-Snaith analogs.

 \enddocument